\documentclass[12pt,a4paper,twoside]{article}

\usepackage{amsmath,amsfonts,amsthm,amsopn,color,amssymb}
\usepackage{palatino}
\usepackage{graphicx}
\newcommand{\C}{\mathbb{C}}
\newcommand{\R}{\mathbb{R}}

\newcommand{\RT}{{\mathbb{R}^3}}

\renewcommand{\le}{\leqslant}
\renewcommand{\ge}{\geqslant}

\newcommand{\vfi}{\varphi}

\renewcommand{\l }{\lambda}
\newcommand{\n }{\nabla }
\newcommand{\s }{\sigma }

\renewcommand{\O}{\Omega}

\newcommand{\Ne}{\mathcal{N}}

\renewcommand{\C}{\mathbb{C}}
\renewcommand{\o}{\omega}

\newcommand{\D }{{\mathcal D}^{1,2}(\RT)}

\renewcommand{\H}{H^{1}(\RT)}

\newcommand{\irt }{\int_{\RT}}

\def\bbm[#1]{\mbox{\boldmath $#1$}}

\newtheorem{theorem}{Theorem}[section]
\newtheorem{lemma}[theorem]{Lemma}

\renewenvironment{proof}{\noindent{\textbf{Proof\quad}}}{$\hfill\square$\vspace{0.2 cm}\\}

\title{{\bf Ground state solutions for the nonlinear Klein-Gordon-Maxwell equations}\footnote{Sponsored by Research Funds PRIN07 ``Metodi Variazionali e topologici nello studio di fenomeni nonlineari"}}

\author{A. Azzollini \footnote{Dipartimento di Matematica, Universit\`a degli
Studi di Bari,  Via E. Orabona 4, I-70125 Bari, Italy, e-mail: {\tt azzollini@dm.uniba.it}}
 \; \& \;
A. Pomponio\footnote{Dipartimento di Matematica, Politecnico di Bari, Via Amendola 126/B, I-70126 Bari, Italy, e-mail:
{\tt a.pomponio@poliba.it}}}

\date{}

\begin{document}

\maketitle

\begin{abstract}
In this paper we prove the existence of a ground state solution for the nonlinear Klein-Gordon-Maxwell equations in the electrostatic case.
\end{abstract}

\section{Introduction}

In this paper we are interested in studying the following
nonlinear Klein-Gordon-Maxwell equations
\begin{equation} \label{eqq}\tag{${\cal KGM}$}
\left\{
\begin{array}{l}
\square u
+\left[  \left\vert \nabla S-e\mathbf{A}\right\vert^{2}
-\left( \frac{\partial S} {\partial t}+e\varphi\right)^{2}+m_0^{2}\right]\,u
-\left\vert u\right\vert ^{p-1}u=0;
\\
\frac{\partial}{\partial t}
\left[  \left(  \frac{\partial S}{\partial t}+e\varphi\right)u^{2}\right]
-\nabla\cdot\left[  \left(  \nabla S-e\mathbf{A}\right)u^{2}\right]  =0;
\\
\nabla\cdot\left(\frac{\partial \mathbf{A}}{\partial t}
+\nabla\varphi\right)
=e\left(  \frac{\partial S}{\partial t}+e\varphi\right)  u^{2};
\\
\nabla\times\left(  \nabla\times\mathbf{A}\right)
+\frac{\partial}{\partial t}\left(  \frac{\partial \mathbf{A}}{\partial t}
+\nabla\varphi\right)
=e\left(  \nabla S-e\mathbf{A}\right) u^{2}.
\end{array}
\right.
\end{equation}
where $e,m_0>0,$ $1< p<5,$ $u(x,t)\in\R,$ $S(x,t)\in\R,$
$(\phi(x,t),\mathbf{A}(x,t))\in\R\times\R^3$. This system arises
in a very interesting physical context: in fact, it provides a
``dualistic model" for the description of the interaction between
a charged relativistic particle of matter and the electromagnetic
field that it generates. According to such a model, the matter
particle is a solitary wave $u(x,t) e^{i S(x,t)}$ which is solution
of a nonlinear field equation, and the interaction with the
electromagnetic field described by the gauge potentials
$(\phi,\mathbf{A})$ is obtained by coupling the field equation
with the Maxwell equations (see \cite{BF2}).

By the invariance of the system with respect to the group of transformations of
Poincar\'e, in order to find a solitary wave it is sufficient to
look for a ``standing wave'' $u(x)e^{i\o t}$ (here $\o\in\R$) and
to make it travel by means of a Lorentz transformation. The
existence of standing waves for \eqref{eqq}, which has been proved
recently by Benci \& Fortunato in \cite{BF2} and D'Aprile \&
Mugnai in \cite{DM}, is a consequence of the nonlinear structure
of the system. In fact, it is well known that in general wave
equations do not possess solitary wave solutions. A typical
example is the Klein-Gordon equation
\[
\square \psi+m^2\psi=0,\quad m\neq 0,\quad\psi(x,t)\in\C,
\]
whose solutions have a spreading behavior which is time dependent
(see \cite{Kic}).

The characteristic of the solitary waves of preserving their
energy density as a localized packet which travels as time goes
on, makes the solitary waves behavior similar to that of the
particle. Differently from the classical model, where the particle
is represented as a dimensionless point, here the particle is
endowed with space extension and has finite energy. This fact
allows us to avoid the well known problem of the {\it divergence
of the energy} which, in the theory of special relativity, brings
to the impossibility of describing the dynamics of the particle
(in fact the inertial mass is infinite: see for example \cite{J},
\cite{LL} and \cite{T}). This is the reason why the solitary waves
appear in several mathematical physics contexts, such as classical
and quantum field theory, nonlinear optics, fluid mechanics,
plasma physics (see e.g. \cite{D.E.G.M.}, \cite{F}, \cite{Kic},
\cite{Ra}, \cite{Wi}).

Finally, it is quite a remarkable fact that, since \eqref{eqq} is
invariant with respect to the the Poincar\'{e} group of
transformations, the model described by \eqref{eqq} turns out to
be consistent with the basic principles of special relativity
theory (see \cite{B} and \cite{BF1}). As a consequence, the
solitary waves  experience well known relativistic phenomenona
such as length contraction, time dilatation and the equivalence
between mass and energy.

In this paper, we are interested in looking for ground state
solutions of the electrostatic \eqref{eqq}, namely for solutions
which minimizes the action among all the solutions. The interest
in ground states, which has been emphasized in many papers such as
the celebrated works of Coleman, Glaser \& Martin \cite{CGM}
and of Berestycki \& Lions \cite{BL1}, is justified by the fact
that they in general exhibit some type of stability. From a
physical point of view, the stability of a
standing wave is a crucial point to establish the existence of soliton-like solutions.

A first work in this direction is the recent paper of Long (see
\cite{Lo}), where the stability properties of the solutions of
\eqref{eqq} have been investigated, for $e$ sufficiently small.

Consider the system \eqref{eqq} in its electrostatic form, namely
set $\mathbf{A}=0,$ $\phi(x,t)=\phi(x),$ $u(x,t)= u(x)$ and
$S(x,t)=\o t:$
\begin{equation}    \label{eq}
\left\{
\begin{array}{ll}
-\Delta u +[m^2_0-(\o+e\phi)^2 ]u-|u|^{p-1}u=0 & \hbox{in }\RT
\\
-\Delta \phi+e^2 u^2 \phi=-e\o u^2 & \hbox{in }\RT.
\end{array}
\right.
\end{equation}

Solutions of \eqref{eq}, $(u,\phi)\in \H\times\D$, are critical points of the functional ${\cal S}:\H\times\D \to \R$ defined as
\[
{\cal S}(u,\phi)=\frac12 \irt |\n u|^2 -|\n \phi|^2+[m^2_0-(\o+e\phi)^2 ]u^2
-\frac{1}{p+1}\irt |u|^{p+1}.
\]

We are interested in finding ``ground state'' solutions of
\eqref{eq}, that is a solution $(u_0,\phi_0)\in \H\times \D$ which
minimizes the functional ${\cal S}$ among all the non-trivial
solutions of \eqref{eq}, namely ${\cal S}(u_0,\phi_0) \le {\cal
S}(u,\phi)$, for any $(u,\phi)\neq (0,0)$ solution of \eqref{eq}.

The main result we provide in this paper is the following
\begin{theorem}\label{main}
The problem \eqref{eq} admits a ground state solution if
\begin{itemize}
    \item $3\le p<5$ and $m_0>\o$;
    \item $1<p<3$ and $m_0\sqrt{p-1}>\o\sqrt{5-p}$.
\end{itemize}
\end{theorem}

\vspace{0.5cm}
\begin{center}
{\bf NOTATION}
\end{center}

\begin{itemize}
\item For any $1\le s< +\infty$, $L^s(\RT)$ is the usual Lebesgue space endowed with the norm
\[
\|u\|_s^s:=\irt |u|^s;
\]
\item $\H$ is the usual Sobolev space endowed with the norm
\[
\|u\|^2:=\irt |\n u|^2+ u^2;
\]
\item $\D$ is completion of $C_0^\infty(\RT)$ (the compactly
supported functions in $C^\infty(\RT)$) with respect to the norm
\[
\|u\|_{\D}^2:=\irt |\n u|^2;
\]
\item for any $r>0,$ $x\in\RT$ and $A\subset \RT$
\begin{align*}
B_r(x) &:=\{y\in\RT\mid |y-x|\le r\},\\
B_r    &:=\{y\in\RT\mid |y|\le r\},\\
A^c    &:= \RT\setminus A.
\end{align*}
\end{itemize}

\section{Preliminary lemmas}

The first difficulty in dealing with the functional ${\cal S}$ is
that it is strongly indefinite, namely it is unbounded both from
below and from above on infinite dimensional subspaces. To avoid
this indefiniteness, we will use the reduction method.

We need the following:
\begin{lemma}\label{le:phipsi}
For any $u\in \H$, there exists a unique $\phi=\phi_u\in \D$ which satisfies
\[
-\Delta \phi+e^2 u^2 \phi=-e\o u^2 \quad \hbox{in }\RT.
\]
Moreover, the map $\Phi:u\in\H\mapsto\phi_u\in\D$ is continuously
differentiable, and on the set $\{x\in \RT\mid u(x)\neq 0\}$,
\begin{equation}\label{eq:phi}
-\frac{\o}{e}\le \phi_u\le 0.
\end{equation}
\end{lemma}

\begin{proof}
The proof can be found in \cite{BF2,DM2}.
\end{proof}

\begin{lemma}\label{le:phipsi2}
Let $u\in\H$ and set $\psi_u =(\Phi'[u])[u]/2\in\D.$\\
Then:
\begin{itemize}
\item $\psi_u$ is a solution to the integral equation
\begin{equation}\label{eq:psiu}
\irt e\o \psi_u u^2= \irt e(\o + e\phi_u)\phi_u u^2;
\end{equation}
\item it results that
\begin{equation}\label{eq:psi}
\psi_u\le 0.
\end{equation}
\end{itemize}
\end{lemma}

\begin{proof}
The proof is a consequence of the fact that $\psi_u$ satisfies
\[
-\Delta \psi_u+ e^2 u^2 \psi_u= - e (\o+e\phi_u)u^2,
\]
as we know by \cite{DM2}.
\end{proof}

Set $\O=m^2_0-\o^2$ and define $I:\H\to \R$ as
\[
I(u)=\frac12 \irt |\n u|^2 +\O u^2 -e\o \phi_u u^2
-\frac{1}{p+1}\irt |u|^{p+1},
\]

The functional $I$ is obtained from ${\cal S}$ by the reduction
method, as in \cite{BF2}. As one can see, it does not present
anymore the strong indefiniteness, and it is strictly connected
with our problem, since $(u,\phi)\in \H\times \D$ is a solution of
\eqref{eq} if and only if $u$ is a critical point of $I$ and
$\phi=\phi_u$.

We will look for a minimizer of the functional $I$ restricted to
the its Nehari manifold, namely
\[
\Ne=\{u\in \H\setminus \{0\}\mid G(u)=0\},
\]
where
\[
G(u)=\langle I'(u),u\rangle= \irt |\n u|^2 +\O u^2 -2e\o \phi_u u^2 -e^2\phi_u^2 u^2-\irt |u|^{p+1}.
\]

In the following Lemmas, we point out some properties related with
the Nehari manifold

\begin{lemma}\label{le:p+1>0}
There exists a positive constant $C$ such that $\|u\|_{p+1}\ge C$, for all $u\in \Ne$.
\end{lemma}

\begin{proof}
By \eqref{eq:phi}, we infer
\[
-e\irt (2\o +e\phi_u)\phi_u u^2\ge 0.
\]
Therefore, by the definition of the Nehari manifold, we get
\[
\|u\|_{p+1}^2\le  C \irt |\n u|^2 +\O u^2 \le C \|u\|_{p+1}^{p+1}.
\]
\end{proof}

\begin{lemma}\label{le:pos}
There exists a positive constant $C>0$, such that $I(u)\ge C$, for any $u\in \Ne$.
\end{lemma}

\begin{proof}
For any $u\in \Ne$, we have
\begin{equation}\label{eq:IN}
I(u) = \frac {p-1}{2(p+1)}\irt |\n u|^2+\O u^2
-\frac {p-3}{2(p+1)}\irt e\o\phi_u u^2+\frac 1 {p+1} \irt e^2\phi_u^2 u^2.
\end{equation}
We have to distinguish two cases. If $3\le p<5$, then, by \eqref{eq:phi}, each term in \eqref{eq:IN} is positive and the conclusion follows by Lemma \ref{le:p+1>0}, supposing $m_0>\o$.
\\
Instead, in the case $1<p<3$, by \eqref{eq:phi} we have
\begin{align*}
I(u)&\ge \frac {p-1}{2(p+1)}\irt |\n u|^2 +\O u^2
+\frac {p-3}{2(p+1)}\irt \o^2 u^2
\\
&\ge \frac {p-1}{2(p+1)}\irt |\n u|^2
+\frac {1}{2(p+1)}\irt [(p-1)m_0^2-2\o^2] u^2.
\end{align*}
Assuming that $m_0\sqrt{p-1}>\o\sqrt{5-p}$, we conclude also in this case.
\end{proof}

\begin{lemma}\label{le:N}
$\Ne$ is a $C^1$ manifold.
\end{lemma}

\begin{proof}
For all $u\in \H$, we have
\[
G(u)=2I(u)+ \irt \frac{1-p}{p+1}|u|^{p+1}-\irt e \o \phi_u u^2 -\irt e^2 \phi_u^2 u^2.
\]
Let us prove that there exists $C>0$ such that $\langle G'(u),u\rangle \le -C$, for all $u\in \Ne$.
\\
If $u\in \Ne$, by \eqref{eq:psiu}
\begin{align*}
\langle G'(u),u\rangle
&= \irt (1-p)|u|^{p+1}-\irt 4 e \phi_u u^2(\o+e \phi_u +e\psi_u)
\\
&=(1-p)\irt |\n u|^2 +\O u^2
\\
&\quad-\irt e\phi_u u^2 [(1-p)(2\o+e\phi_u)+4(\o+e \phi_u +e\psi_u)].
\end{align*}
We have to distinguish two cases. If $3\le p<5$, since $m_0>\o$,
by Lemma \ref{le:p+1>0} and \eqref{eq:phi}, we need only to show that
\begin{equation*}
(1-p)(2\o+e\phi_u)+4(\o+e \phi_u +e\psi_u)\le 0.
\end{equation*}
Indeed, since $\phi_u,\psi_u\le0$, we have
\begin{equation*}
(1-p)(2\o+e\phi_u)+4(\o+e \phi_u +e\psi_u)
=2(3-p)\o
+(5-p)e\phi_u
+4e\psi_u \le 0.
\end{equation*}
In the case $1<p<3$, instead, by \eqref{eq:phi}, we have
\begin{align*}
\langle G'(u),u\rangle
&\le (1-p)\irt |\n u|^2 +\O u^2
- 2(3-p)\irt e\o\phi_u u^2
\\
&\quad -(5-p)\irt e^2 \phi_u^2 u^2
-4\irt e^2\phi_u \psi_u u^2
\\
&\le (1-p)\irt |\n u|^2
+\irt [(1-p)m_0^2+(5-p)\o^2]u^2.
\end{align*}
We get the same conclusion with the additional assumption $m_0\sqrt{p-1}>\o\sqrt{5-p}$.
\end{proof}

According to the definition of \cite{L1}, we say that a sequence $(v_n)_n$ vanishes if, for all $r>0$
\[
\lim_n \sup_{\xi \in \RT} \int_{B_r(\xi)}v_n^2=0.
\]
\begin{lemma}\label{le:nonvan}
Any bounded sequence $(v_n)_n\subset \Ne$ does not vanish.
\end{lemma}

\begin{proof}
Suppose by contradiction that $(v_n)_n$ vanishes, i.e. there
exists $\bar r>0$ such that
\[
\lim_n \sup_{\xi \in \RT} \int_{B_{\bar r}(\xi)}v_n^2=0.
\]
Then, by \cite[Lemma 1.1]{L2}, we infer that $v_n\to 0$ in $L^s(\RT)$, for any $2<s<6$, contradicting Lemma \ref{le:p+1>0}.
\end{proof}

The map $\Phi$ is continuous for the weak topology in the sense of
the following lemma
    \begin{lemma}\label{le:weakcon}
        If $u_n\rightharpoonup u_0$ in $\H$ then, up to subsequences, $\phi_{u_n}\rightharpoonup\phi_{u_0}$ in
        $\D$. As a consequence
        $I'(u_n)\to I'(u_0)$ in the sense of distributions.
    \end{lemma}
    \begin{proof}
        Let $(u_n)_n$ and $u_0$ be in $\H,$ and assume that
        $u_n\rightharpoonup u_0$ in $\H.$ As a consequence
            \begin{align}
                u_n &\rightharpoonup u_0, \hbox{ in } L^s(\RT),\; 2\le s\le 6,\label{eq:w}\\
                u_n &\to u_0, \hbox{ in }L^s_{loc}(\RT),\; 1\le s<6.\label{eq:loc}
            \end{align}
    We denote by $\phi_n$ the function $\phi_{u_n}.$
        By the second of
        \eqref{eq} we have that for any $n\ge 1$
            \begin{align*}
              \irt|\n \phi_n|^2
              &= - e^2 \irt u_n^2 \phi_n^2 - e \irt\o u_n^2\phi_n\\
              &\le - e \irt\o u_n^2\phi_n\le C \|u_n\|_{12/5}^2\|\n
              \phi_n\|_2,
            \end{align*}
        and then we deduce that $(\phi_n)_n$ is bounded in $\D$.
        We can assume that there exists $\phi_0\in\D$ such that $\phi_n\rightharpoonup
        \phi_0$ in $\D$ and, as a consequence,
            \begin{align}
                \phi_n &\rightharpoonup \phi_0, \hbox{ in } L^6(\RT), \label{eq:w2}\\
                \phi_n &\to \phi_0, \hbox{ in }L^s_{loc}(\RT),\; 1\le s<6.\label{eq:loc2}
            \end{align}
        If we show that $\phi_0=\phi_{u_0}$ we have
        concluded. By the uniqueness of the solution of the second
        equation in \eqref{eq}, we are reduced to prove that
            \begin{equation*}
                -\Delta \phi_0+e^2 u_0^2 \phi_0=-e\o u_0^2
            \end{equation*}
        in the sense of distributions.\\
        So, let
        $\varphi\in C_0^\infty(\RT)$ a test function. Since
            $$-\Delta \phi_n+e^2 u_n^2 \phi_n=-e\o u_n^2$$
        it is sufficient to show that the following three hold
            \begin{align}
                \irt (\n\phi_n|\n\varphi)&\to\irt
                (\n\phi_0|\n\varphi)\nonumber\\
                \irt u_n^2\phi_n\varphi&\to\irt
                u_0^2\phi_0\varphi\label{eq:sec}\\
                \irt u_n^2\varphi&\to\irt
                u_0^2\varphi.\nonumber
            \end{align}
        The first is a trivial application of the definition of
        weak convergence, whereas the third is a consequence of \eqref{eq:loc}. As regards the
        second,
        observe that
            \begin{align*}
                \irt (u_n^2\phi_n-
                u_0^2\phi_0)\varphi &= \irt (u_n^2-
                u_0^2)\phi_n\varphi + \irt (\phi_n-\phi_0)
                u_0^2\varphi\\
                &\le C \|\n \phi_n\|_2\left(\irt |u_n^2-u_0^2|^{\frac 6
                5}|\varphi|^{\frac 6 5}\right)^{\frac 5 6}+\irt (\phi_n-\phi_0)
                u_0^2\varphi
            \end{align*}
        and then \eqref{eq:sec} follows by the boundedness of
        $(\phi_n)_n$, \eqref{eq:loc} and \eqref{eq:loc2}.\\
        Now we pass to prove the second part of the Lemma. Let $\varphi$ be a test function.\\ We
        compute:
            \begin{align*}
                \langle I'(u_n), \varphi\rangle & = \irt (\n u_n|\n\varphi) +\O u_n\varphi -2e\o \phi_n
                u_n\varphi
                -e^2\phi_n^2 u_n\varphi- |u_n|^{p-1}u_n\varphi\\
                \langle I'(u_0), \varphi\rangle & = \irt (\n u_0|\n\varphi) +\O u_0\varphi -2e\o
                \phi_0
                u_0\varphi
                -e^2\phi_0^2 u_0\varphi- |u_0|^{p-1}u_0\varphi.
            \end{align*}
        Now observe that
            \begin{align*}
                \irt (\phi_n u_n - \phi_0 u_0)\varphi&=
                \irt \phi_n(u_n-u_0)\varphi + \irt (\phi_n-\phi_0)u_0\varphi\\
                &\le C\|\n\phi_n\|_2\left(\irt |u_n-u_0|^{\frac 6
                5}|\varphi|^{\frac 6 5}\right)^{\frac 5 6} + \irt
                (\phi_n-\phi_0)u_0\varphi\\
                &= o_n(1)
            \end{align*}
        by the boundedness of
        $(\phi_n)_n$, \eqref{eq:loc} and \eqref{eq:loc2}. Moreover,
            \begin{align*}
                \irt (\phi_n^2 u_n - \phi^2_0 u_0)\varphi&=
                \irt \phi^2_n(u_n-u_0)\varphi + \irt (\phi^2_n-\phi^2_0)u_0\varphi\\
                &\le C\|\n\phi_n\|^2_2\left(\irt |u_n-u_0|^{\frac
                3 2}|\varphi|^{\frac 3 2}\right)^{\frac 2 3} + \irt
                (\phi_n^2-\phi_0^2)u_0\varphi\\
                &= o_n(1)
            \end{align*}
        by the boundedness of $(\phi_n)_n$, \eqref{eq:loc} and
        \eqref{eq:loc2}.
So we have
        \begin{equation*}
\begin{array}{c}
\underbrace{\irt (\n u_n|\n\varphi) +\O u_n\varphi}_{\downarrow} -
\underbrace{\irt2e\o \phi_n u_n\varphi}_{\downarrow} -
\underbrace{\irt e^2\phi_n^2 u_n\varphi}_{\downarrow}-
\underbrace{\irt|u_n|^{p-1}u_n\varphi}_{\downarrow}
\\
\displaystyle \irt (\n u_0|\n\varphi) +\O u_0\varphi
-\irt2e\o\phi_0 u_0\varphi -\irt e^2\phi_0^2 u_0\varphi-\irt
|u_0|^{p-1}u_0\varphi.
\end{array}
\end{equation*}
and then we conclude that $\langle I'(u_n),
\varphi\rangle\to\langle I'(u_0), \varphi\rangle$.
\end{proof}

\section{Proof of Theorem \ref{main}}

Let $\s=\inf_{u\in\Ne} I(u)$. By Lemma \ref{le:pos}, we argue that $\s>0$. Since all the critical points of $I$ are contained in $\Ne$ and since, by Lemma \ref{le:N}, we know that Nehari manifold is a natural constrained for $I$, if there exists $u_0\in \Ne$ such that $I(u_0)=\s$, then $(u_0,\phi_{u_0})$ is a ground state solution for \eqref{eq}.

Let $(u_n)_n\subset \Ne$ such that $I(u_n)\to \s$, as $n\to \infty$. It is easy to see that $(u_n)_n$ is a bounded sequence in $\H$.
By Lemma \ref{le:nonvan}, there exists $C>0$, $\bar r>0$ and a sequence $(\xi_n)_n\subset\RT$ such that
\[
\int_{B_{\bar r}(\xi_n)}u_n^2 \ge C.
\]
Let $v_n=u_n(\cdot + \xi_n)$. By the invariance of translations, $(v_n)_n$ is a bounded sequence contained in $\Ne$ such that
\begin{equation}\label{eq:br}
\int_{B_{\bar r}}v_n^2 \ge C,\quad \hbox{for all }n,
\end{equation}
and, moreover, $I(v_n)\to \s$, as $n\to \infty$.
Up to a subsequence, there exists $v_0\in \H$ such that
\begin{align}
v_n \rightharpoonup v_0,& \hbox{ weakly in }\H,\nonumber
\\
v_n \to v_0,& \hbox{ in }L^s_{loc}(\RT),\; 1\le s<6,\nonumber
\\
v_n \to v_0,& \hbox{ a.e. in }\RT,\label{eq:q.o.}
\end{align}
Denote $\phi_n\equiv \phi_{v_n}$ and $\phi_0\equiv\phi_{v_0}$.\\
By Lemma \ref{le:weakcon}, we know that
$\phi_n\rightharpoonup\phi_0$ in $\D$, and, as a consequence,
            \begin{align}
                \phi_n &\to \phi_0, \hbox{ in }L^s_{loc}(\RT),\; 1\le
                s<6\nonumber\\
                \phi_n & \to \phi_0, \hbox{ a.e. in
                }\RT.\label{eq:qo}
            \end{align}
By \cite{Wi}, without lost of generality, we can assume that $(v_n)_n$ is a Palais-Smale
sequence for the functional $I_{|\Ne}$, in particular,
\begin{align}
I(v_n) \to \s, &\hbox{ as }n\to \infty, \nonumber
\\
(I_{|\Ne})'(v_n)\to 0, &\hbox{ as }n\to \infty. \label{eq:psn}
\end{align}
By \eqref{eq:psn}, being $(v_n)_n$ bounded in $\H$, for suitable Lagrange multipliers $\l_n$, we get
\[
o_n(1)=\langle(I_{|\Ne})'(v_n),v_n\rangle
=\langle I'(v_n),v_n\rangle
+\l_n \langle G'(v_n),v_n\rangle
=\l_n \langle G'(v_n),v_n\rangle.
\]
By Lemma \ref{le:N}, we infer that $\l_n=o_n(1)$ and, by
\eqref{eq:psn},
    \begin{equation}\label{eq:psnn}
        I'(v_n)\to 0, \hbox{ as }n\to \infty.
    \end{equation}

By \eqref{eq:br}, we infer that $v_0\neq 0$ (and hence also
$\phi_0\neq 0$). Moreover, by Lemma \ref{le:weakcon} and
\eqref{eq:psnn}, we can conclude that $I'(v_0)=0.$ It remains to
prove that $I(v_0)=\s$. Observe that, since $(v_n)_n$ is in $\Ne$,
we have
    \begin{equation*}
        I(v_n) = \frac {p-1}{2(p+1)}\irt |\n v_n|^2
        +\O v_n^2 -
        \frac {p-3}{2(p+1)}\irt
        e\o\phi_nv_n^2+\frac 1 {p+1} \irt e^2\phi_n^2 v_n^2.
    \end{equation*}
We have to distinguish two cases. If $p\ge 3$, since $\phi_n\le 0$, by the weak lower
semicontinuity of the $H^1-$norm, \eqref{eq:q.o.}, \eqref{eq:qo}
and the Lemma of Fatou, we conclude that $I(v_0)=\s.$ This implies
that $(v_0,\phi_0)$ is a ground state solution. If $1<p<3$, by \eqref{eq:phi} and requiring that $m_0\sqrt{p-1}>\o\sqrt{5-p}$, it is easy to see that
\[
\frac {p-1}{2(p+1)}\O v_n^2 - \frac {p-3}{2(p+1)} e\o\phi_n v_n^2 \ge 0, \qquad \hbox{a.e. in }\RT,
\]
and we conclude as before.


\begin{thebibliography}{99}

\bibitem{B}
V. Benci, {\it Some remarks on the theory of relativity and the
naive realism} In: Cerrai, P. et al. (eds.): The applications of
mathematics to the sciences of nature. Kluwer, Dordrecht (2002),
39--78.

\bibitem{BF}
V. Benci, D. Fortunato, {\it An eigenvalue problem for the
Schr\"odinger-Maxwell equations},  Topol. Methods Nonlinear Anal.,
{\bf 11} (1998), 283--293.

\bibitem{BF1}
V. Benci, D. Fortunato, {\it Some remarks on the semilinear wave
equations}, Progress in Nonlinear Differential equations and Their
Applications {\bf 54} Birkhauser, Basel (2003), 141--162.

\bibitem {BF2}
V. Benci, D. Fortunato, {\it Solitary waves of the nonlinear
Klein-Gordon field equation coupled with the Maxwell equations,}
Rev. Math. Phys. {\bf 14} (2002), 409--420.

\bibitem{BFMP}
V. Benci, D. Fortunato, A. Masiello, L. Pisani,
{\it Solitons and the electromagnetic field},  Math. Z., {\bf 232}, (1999), 73--102.


\bibitem{BL1}
H. Berestycki, P.L. Lions, {\it Nonlinear scalar field equations,
I - Existence of a ground state}, Arch. Rational Mech. Anal., {\bf
82}, (1983), 313--345.

\bibitem{BL2}
H. Berestycki, P.L. Lions, {\it Nonlinear scalar field fquations,
II - Existence of infinitely many solutions}, Arch. Rational Mech.
Anal., {\bf 82}, (1983), 347--375.

\bibitem{C}
D. Cassani, {\it Existence and non-existence of solitary waves for
the critical Klein-Gordon equation coupled with Maxwell's
equations}, Nonlinar Anal., {\bf 58} (7-8) (2004), 733--747.



\bibitem {CGM}
S. Coleman, V. Glaser, A.Martin, {\it Action minima amoung
solutions to a class of euclidean scalar field equations}, Commun.
math. Phys., {\bf 58}, (1978), 211--221.

\bibitem {D.E.G.M.}
K. Dodd, D. Eilbeck, J.D. Gibbon, H.C. Morris \textit{\emph{
Solitons and Nonlinear Wave Equations, }} \emph{Academic Press,
London, New York}, (1982).

\bibitem {DM}
T. D'Aprile, D. Mugnai, {\it Solitary waves for nonlinear
Klein-Gordon-Maxwell and Schr\"{o}dinger-Maxwell equations},
Proc. Roy. Soc. Edinburgh Sect. A, {\bf 134}, (2004), 893--906.

\bibitem{DM2}
T. D'Aprile, D. Mugnai, {\it Non-existence results for the coupled Klein-Gordon-Maxwell equations},
Adv. Nonlinear Stud., {\bf 4}, (2004), 307--322.

\bibitem {F}
B. Felsager, \textit{\emph{Geometry, Particle and fields,}}
\emph{Odense University press}, (1981).

\bibitem {J}
J.D. Jackson, {\it\emph{Classical electrodynamics, }} \emph{John
Wiley \& Sons., New York, London}, (1962).

\bibitem{Kic}
S. Kichenassamy,\:\textit{\emph{Nonlinear wave equations}},
\emph{Dekker, NY}, 1996.

\bibitem {LL}
L. Landau, E. Lifchitz, {\it\emph{Th\'{e}orie des Champs}},
\emph{Editions Mir, Moscou,} (1970).

\bibitem{L1}
P.L. Lions,
\textit{The concentration-compactness principle in the calculus of variation.
The locally compact case. Part I},
Ann. Inst. Henri Poincar\'e, Anal. Non Lin\'eaire, {\bf 1}, (1984), 109--145.


\bibitem{L2}
P.L. Lions,
\textit{The concentration-compactness principle in the calculus of variation.
The locally compact case. Part II},
Ann. Inst. Henri Poincar\'e, Anal. Non Lin\'eaire, {\bf 1}, (1984), 223--283.

\bibitem{Lo}
E. Long, {\it Existence and stability of solitary waves in
non-linear Klein-Gordon-Maxwell equations}, Rev. Math. Phys., {\bf
18} (2006), 747--779.

\bibitem {Ra}
R. Rajaraman, \textit{\emph{Solitons and instantons,}} \emph{North
Holland, Amsterdam, Oxord, New York, Tokio,} (1988).

\bibitem {T}
W. Thirring, {\it\emph{Classical Mathematical Physics,}}
\emph{Springer, New York, Vienna,} (1997).

\bibitem{Wi}
M. Willem,
Minimax Theorems.
Progress in Nonlinear Differential Equations and their Applications, 24. Birkh\"auser Boston, Inc., Boston, MA, 1996.

\end{thebibliography}
\end{document}